\documentclass[12pt]{article}
\usepackage{amssymb}
\usepackage{mathrsfs}
\usepackage{pifont}
\usepackage{tipa}
\usepackage{amsmath}
\usepackage{amssymb,amsmath}
\usepackage{amsfonts}
\catcode`@=11
 \long\def\@makefntext#1{\noindent #1}

\begin{document}
\title{\large{\textbf{Unitary modules for the
twisted Heisenberg-Virasoro algebra}} }
\author{Xiufu Zhang $^{1*}$, Shaobin Tan $^2$, Haifeng Lian $^3$\\
{\scriptsize 1. School of Mathematical Sciences, Xuzhou Normal
University, Xuzhou 221116,  China}\\
{\scriptsize 2. School of Mathematical Sciences, Xiamen University,
Xiamen 361005, China}\\
{\scriptsize 3. School of  Computer and Information, Fujian
Agriculture and Forestry University, Fuzhou 350002, China}}
\date{}
\maketitle \footnotetext{\footnotesize Supported by the National
Natural Science Foundation of China (No. 10931006).}
\footnotetext{\footnotesize *Corresponding author. Email:
xfzhang@xznu.edu.cn;\quad tans@xmu.edu.cn;\quad hlian@fjau.edu.cn}

\numberwithin{equation}{section}

\begin{abstract}
In this paper, the conjugate-linear anti-involutions and the
unitary irreducible modules of the intermediate series over the
twisted Heisenberg-Virasoro algebra are classified respectively.
We prove that any unitary irreducible module of the intermediate
series over the twisted Heisenberg-Virasoro algebra is of the form
$\mathcal{A}_{a,b,c}$ for $a\in \mathbb{R}, b\in
\frac{1}{2}+\sqrt{-1}\mathbb{R}, c\in \mathbb{C}.$
\vspace{2mm} \vspace{2mm}
\end{abstract}

\vskip 3mm \noindent{\section{{Introduction}}}

The twisted Heisenberg-Virasoro algebra $\mathfrak{L}$ is defined to
be a Lie algebra with $\mathbb{C}$-basis $\{L_m,I_m,C_L,C_I,C_{LI}
\mid m\in \mathbb{Z}\}$ subject to the following Lie brackets:
\begin{align*}
&[L_m,L_n]=(n-m)L_{m+n}+\delta_{m+n,0}\frac{m^{3}-m}{12}C_L,\\
&[I_m,I_n]=n\delta_{n+m,0}C_I,\\
&[L_m,I_n]=nI_{m+n}+\delta_{m+n,0}(m^{2}-m)C_{LI},\\
&[\mathfrak{L},C_L]=[\mathfrak{L},C_I]=[\mathfrak{L},C_{LI}]=0.
\end{align*}
This Lie algebra was first introduced by Arbarello et al. in Ref.
[1]. It is the universal central extension of the Lie algebra of
differential operators on a circle of order at most one:
$$L_{HV}=\{f(t)\frac{d}{dt}+g(t)|f,g\in\mathbb{C}[t,t^{-1}]\}.$$
By the definition, both  the Heisenberg algebra and the Virasoro
algebra are subalgebras of the twisted Heisenberg-Virasoro algebra
$\mathfrak{L}$. It is known that the twisted Heisenberg-Virasoro
algebra has some important applications in the representation
theory of the toroidal Lie algebra which is the prime example of
the generalizations of Kac-Moody algebra (see [6]). Moreover,
$\mathfrak{L}$ has some relations with the well known $N=2$
Neveu-Schwarz superalgebra. In fact, the even part of the $N=2$
Neveu-Schwarz superalgebra is essentially the twisted
Heisenberg-Virasoro algebra (see [11]).

The representation theories of the Virasoro algebra and its
related algebras play crucial roles in many areas of Mathematics
and Physics and have been well developed.  In particular, the
unitary representations are significant. The unitary highest
weight representations of Virasoro algebras are determined in
[4,5,8-10]. The complete classification of the unitary
Harish-Chandra modules over Virasoro algebra is shown in [3]. The
irreducible Harish-Chandra modules of Virasoro algebras are
classified in [15]. Recently, the representation theories over the
twisted Heisenberg-Virasoro algebra were studied by several
authors.  When the central element of the Heisenberg subalgebra
acts trivially, the highest weight modules for twisted
Heisenberg-Virasoro algebra are studied in [1] and [2]. It is
proved in [1] that the unitary highest weight modules for twisted
Heisenberg-Virasoro algebra is just the unitary highest weight
modules for Virasoro algebra when the central charge of the
Heisenberg subalgebra is zero (Theorem 6.6 (I) in [1]).The
irreducible Harish-Chandra modules, the module of the intermediate
series, the irreducible weight modules with a finite dimensional
weight space and the Whittaker modules over the twisted
Heisenberg-Virasoro algebra are are also studied (see
[1,2,7,12-14,16]).

The goal of the present paper is to study the conjugate-linear
anti-involutions and the unitary irreducible modules of the
intermediate series over the twisted Heisenberg-Virasoro algebra.

The paper is organized as follows. In Section 2, we study the
conjugate-linear anti-involution of the twisted
Heisenberg-Virasoro algebra $\mathfrak{L}$.  In Section 3, the
unitary irreducible Harish-Chandra modules of the intermediate
series are classified.

Throughout this paper we make a convention that the weight modules
over twisted Heisenberg-Virasoro algebra and Virasoro algebra are
all with finite dimensional weight spaces, i.e., the
Harish-Chandra modules. The symbols $\mathbb{C},
\mathbb{R},\mathbb{R}^+, \mathbb{Z}$ and $S^1$ represent for the
complex field, real number field, the set of positive real number,
the set of integers and the set of complex number of modulus one
respectively.

\vskip 3mm\noindent{\section{Conjugate-linear anti-involutions of
$\mathfrak{L}$}}

It is easy to see the following facts about $\mathfrak{L}:$

(1)
$\mathfrak{C}:=\mathbb{C}C_L\oplus\mathbb{C}C_I\oplus\mathbb{C}C_{LI}\oplus\mathbb{C}I_0$
is the center of $\mathfrak{L}.$

(2) If $x\in\mathfrak{L}$ acts semisimply on $\mathfrak{L}$ by the
adjoint action, then $$x\in
\mathfrak{h}:=\mathrm{span}_{\mathbb{C}}\{L_0,
I_0,C_L,C_I,C_{LI}\},$$ the Cartan subalgebra of $\mathfrak{L}$.

(3) $\mathfrak{L}$ has a weight space decomposition according to the
Cartan subalgebra $\mathfrak{h}:$
$$\mathfrak{L}=\bigoplus_{n\in \mathbb{Z}}
\mathfrak{L}_n,$$ where
$\mathfrak{L}_n=\mathrm{span}_{\mathbb{C}}\{L_n, I_n\}$ if $n\neq0$
and $\mathfrak{L}_0=\mathfrak{h}.$

(4) The Heisenberg algebra $H=\oplus_{n\in
\mathbb{Z}}\mathbb{C}I_n\oplus \mathbb{C}C_I$ and the Virasoro
algebra $Vir=\oplus_{n\in \mathbb{Z}}\mathbb{C}L_n\oplus
\mathbb{C}C_L$ are subalgebras of $\mathfrak{L}.$

\vskip 3mm\noindent{\bf{Lemma 2.1.}} $H\oplus \mathbb{C}C_{LI}\oplus
\mathbb{C}C_L$ is the unique maximal ideal of $\mathfrak{L}.$

\vskip 3mm\noindent{\bf{Proof.}} The proof is similar as that for
Lemma 2.2 in [17], we omit the details.\hfill$\Box$

\vskip 3mm\noindent{\bf{Definition 2.2.}} Let $\mathfrak{g}$ be a
Lie algebra and $\theta$ be a conjugate-linear anti-involution of
$\mathfrak{g}$, i.e. $\theta$ is a map
$\mathfrak{g}\rightarrow\mathfrak{g}$ such that
$$\theta(x+y)=\theta(x)+\theta(y),\ \theta(\alpha x)=\bar{\alpha}\theta(x),$$
$$\theta([x,y])=[\theta(y),\theta(x)],\ \theta^{2}=id$$
for all $x,y\in \mathfrak{g},\ \alpha\in \mathbb{C},$ where $id$ is
the identity map of $\mathfrak{g}.$ A module $V$ of $\mathfrak{g}$
is called unitary if there is a positive definite Hermitian form
$\langle\cdot\ ,\cdot\rangle$ on $V$ such that
$$\langle xu,v\rangle=\langle u,\theta(x)v\rangle$$ for all $u, v\in V, x\in \mathfrak{g}.$

\vskip 3mm\noindent{\bf{Lemma 2.3.}} (Proposition 3.2 in Ref. [3])
Any conjugate-linear anti-involution of $Vir$ is one of the
following types:

(i) $\theta_{\alpha}^{+}(L_n)=\alpha^{n}L_{-n},\
\theta_{\alpha}^{+}(C_L)=C_L,$ for some $\alpha\in \mathbb{R}^{*},$
the set of nonzero real number.

(ii) $\theta_{\alpha}^{-}(L_n)=-\alpha^{n}L_{n},\
\theta_{\alpha}^{-}(C_L)=-C_L,$ for some $\alpha\in S^{1},$ the set
of complex number of modulus one.

\vskip 3mm\noindent{\bf{Lemma 2.4.}} Let $\theta$ be an arbitrary
conjugate-linear anti-involution of $\mathfrak{L}.$ Then

(i) $\theta(\mathfrak{h})=\mathfrak{h}.$

(ii)$$\theta(L_0)=\lambda L_0+\lambda^{'}C,$$
$$\theta(I_m)=\alpha_{m,-\lambda m} I_{-\lambda m}(m\neq0),$$
$$\theta(C_I)=\alpha_{1,-\lambda}\alpha_{-1,\lambda}\lambda C_I,$$
where $\lambda\in \{1,-1\}, \lambda^{'}\in \mathbb{C},$
$C\in\mathfrak{C},$ $\alpha_{m,m}\in S^1,
\overline{\alpha_{m,-m}}\alpha_{-m,m}=1$

(iii) $\theta(H\oplus\mathbb{C}C_{LI})=H\oplus\mathbb{C}C_{LI}.$

\vskip 3mm\noindent{\bf{Proof.}} For (i), we see that
$$
\theta(\mathfrak{C}) =\mathfrak{C}\subset \mathfrak{h}.
$$
On the other hand, identities
$[\theta(L_0),\theta(L_n)]=-n\theta(L_n)$  and
$[\theta(L_0),\theta(I_n)]=-n\theta(I_n)$ imply that $\theta(L_{0})$
acts semisimply on $\mathfrak{L}.$ Thus $\theta(L_{0})\in
\mathfrak{h}.$

For (ii), by (i) and $\theta^2(L_0)=L_0,$ we can assume
\begin{equation} \label{eq:1}
\theta(L_0)=\lambda L_0+\lambda^{'}C
\end{equation}
for some $\lambda\in S^{1}, \lambda^{'}\in \mathbb{C}$ and
$C\in\mathfrak{C}.$ On the other hand, we see that
$\theta(H\oplus\mathbb{C}C_{LI})$ is an ideal of $\mathfrak{L}$
since $H\oplus\mathbb{C}C_{LI}$ is an ideal. Thus we have
$\theta(H\oplus\mathbb{C}C_{LI})\subset H\oplus
\mathbb{C}C_{LI}\oplus \mathbb{C}C_L$ by Lemma 2.1,  and we can
assume
\begin{equation} \label{eq:1}
\theta(I_m)=\sum_{n}\alpha_{m,n}I_n+\beta_mC_{I}+\gamma_mC_{LI}+\zeta_m
C_L\ \mathrm{for}\ m\neq0,
\end{equation}
where $\alpha_{m,n}, \beta_m, \gamma_m, \zeta_m\in\mathbb{C},$ $x\in
H.$ By (2.1), (2.2) and $[\theta(I_m),\theta(L_0)]=m\theta(I_m),$ we
deduce that
\begin{equation} \label{eq:1}
\lambda=\pm 1,
\end{equation}
\begin{equation} \label{eq:1}
\theta(I_m)=\left\{\begin{array}
{r@{\quad\quad}l} \alpha_{m,m}I_m \ \mathrm{if}\  \lambda=-1 \\
\alpha_{m,-m}I_{-m}\ \mathrm{if}\   \lambda=1
\end{array}\right.=\alpha_{m,-\lambda m}I_{-\lambda m}\ \mathrm{for}\
m\neq0,
\end{equation}
where $\alpha_{m,m}\in S^{1},
\overline{\alpha_{m,-m}}\alpha_{-m,m}=1$ since $\theta^2(I_m)=I_m.$
Further more, we have
\begin{equation} \label{eq:1}
\theta(C_I)=[\theta(I_1),\theta(I_{-1})]=\alpha_{1,-\lambda}\alpha_{-1,\lambda}\lambda
C_I.
\end{equation}

For (iii), we see that
\begin{equation} \label{eq:1}
\theta(I_0)=[\theta(L_1),\theta(I_{-1})]=[\theta(L_{1}),\alpha_{-1,\lambda}I_{\lambda}]\in
H\oplus \mathbb{C}C_{LI}.
\end{equation}
\begin{equation} \label{eq:1}
\theta(C_{LI})=\frac{1}{2}([\theta(I_1),\theta(L_{-1})]-\theta(I_0))\in
H\oplus \mathbb{C}C_{LI}.
\end{equation}
From (2.4)-(2.7), we have
$\theta(H\oplus\mathbb{C}C_{LI})=H\oplus\mathbb{C}C_{LI}.$
\hfill$\Box$

\vskip 3mm\noindent{\bf{Proposition 2.5.}} Any conjugate-linear
anti-involution of $\mathfrak{L}$ is one of the following types:
\begin{align*}(\mathrm{i}):\ \ &\theta_{\alpha,\gamma}^{+}(L_n)
=\alpha^{n}L_{-n},\\
&\theta_{\alpha,\gamma}^{+}(C_L)=C_L,\\
&\theta_{\alpha,\gamma}^{+}(I_m)=\alpha^{m}e^{ir}I_{-m}(m\neq0),\\
&\theta_{\alpha,\gamma}^{+}(I_0)=e^{ir}I_0+2e^{ir}C_{LI},\\
&\theta_{\alpha,\gamma}^{+}(C_I)=e^{2ir}C_{I},\\
&\theta_{\alpha,\gamma}^{+}(C_{LI})=-e^{ir}C_{LI}.
\end{align*}
where  $\alpha=\mathbb{R}^*, \gamma\in\mathbb{R}, i=\sqrt{-1}.$
\begin{align*}(\mathrm{ii}):\ \ &\theta^-_{\alpha,\alpha_1,\beta_1,\beta_{-1}}(L_n)
=-\alpha^{n}L_{n}+(\frac{n+1}{2}\alpha^{n-1}\beta_1
-\frac{n-1}{2}\alpha^{n+1}\beta_{-1})I_{n},\\
&\theta^-_{\alpha,\alpha_1,
\beta_1,\beta_{-1}}(C_L)=-C_L-12(\alpha^{-1}\beta_{1}-\alpha\beta_{-1})C_{LI}
+(14\beta_{-1}\beta_{1}-3\alpha^{-2}\beta^{2}_{1}-3\alpha^{2}\beta^{2}_{-1})C_I,\\
&\theta^-_{\alpha,\alpha_1, \beta_1,\beta_{-1}}(I_n)=\alpha_nI_n \ \mathrm{for}\ n\neq0,\\
&\theta^-_{\alpha,\alpha_1, \beta_1,\beta_{-1}}(I_0)=\alpha_{1}\alpha^{-1} I_0-\alpha_{1}\alpha^{-2}\beta_1C_I,\\
&\theta^-_{\alpha,\alpha_1, \beta_1,\beta_{-1}}(C_I)=-\alpha_1^2\alpha^{-2}C_I,\\
&\theta^-_{\alpha,\alpha_1,\beta_1,\beta_{-1}}(C_{LI})=\mu\alpha^{-1}C_{LI}
+\frac{1}{2}(\alpha_{1}\alpha^{-2}\beta_1-\alpha_1\beta_{-1})C_I,\\
\end{align*}
where $\alpha, \alpha_n(n\neq0)\in
S^{1},\beta_1,\beta_{-1}\in\mathbb{C},\ \mathrm{satisfying}\
\alpha_{n}=\alpha_1\alpha^{n-1},
\alpha_1\overline{\beta_1}=\overline{\alpha}\beta_{1},
\alpha\beta_{-1}=\alpha_1\overline{\beta_{-1}},(\overline{\alpha}-\alpha)\beta_{-1,1}=0$
and $\beta_{1,-1}\beta_{-1,1}(1+\overline{\alpha}^2)=0.$

\vskip 3mm\noindent{\bf{Proof.}} Let $\theta$ be any
conjugate-linear anti-involution of $\mathfrak{L}.$ By Lemma 2.4
(iii), we have the induced conjugate-linear anti-involution of
$\mathfrak{L}/(H\oplus \mathbb{C}C_{LI})\simeq Vir:$
$$\bar{\theta}:\mathfrak{L}/(H\oplus\mathbb{C}C_{LI})\rightarrow\mathfrak{L}/(H\oplus\mathbb{C}C_{LI}).$$
Thus, by Lemma 2.3, we see that $\bar{\theta}$ is one of the
following types:

(a) $\bar{\theta}_{\alpha}^{+}(\bar{L_n})=\alpha^{n}\bar{L_{-n}},\
\bar{\theta}_{\alpha}^{+}(\bar{C_L})=\bar{C_L},$ for some $\alpha\in
\mathbb{R}^{*}.$

(b) $\bar{\theta}_{\alpha}^{-}(\bar{L_n})=-\alpha^{n}\bar{L_{n}},\
\bar{\theta}_{\alpha}^{-}(\bar{C_{L}})=-\bar{C_{L}},$ for some
$\alpha\in S^{1}.$

If $\bar{\theta}$ is of type (a), then we can assume
\begin{equation} \label{eq:1}
\theta(L_n)=\alpha^{n}L_{-n}+\sum_{i}\beta_{n,i}I_i+ \gamma_{n}
C_{I}+\zeta_{n} C_{LI},
\end{equation}
where $\alpha\in \mathbb{R}^{*},$ $\beta_{n}, \gamma_{n}, \zeta_n\in
\mathbb{C}.$ By (2.8) and Lemma 2.4 (ii), we have
\begin{equation} \label{eq:1}
\theta(L_0)=L_{0}+\beta_{0,0}I_0+\gamma_0C_I+\zeta_{0} C_{LI}.
\end{equation}
By (2.8)-(2.9) and $[\theta(L_n),\theta(L_0)]=n\theta(L_n),$ we can
deduce easily that $\beta_{n,i}=0$ unless $i=-n,$
$\gamma_{n}=\zeta_n=0$ for all $n\neq0.$ Thus
\begin{equation} \label{eq:1}
\theta(L_n)=\alpha^{n}L_{-n}+\beta_{n,-n}I_{-n}\ \mathrm{for}\
n\neq0.
\end{equation}
By $[\theta(L_{-1}),\theta(L_{1})] =-2\theta(L_0),$ we deduce that
\begin{equation} \label{eq:1}\gamma_0=\frac{1}{2}\beta_{-1,1}\beta_{1,-1},
\end{equation}
\begin{equation} \label{eq:1}
\beta_{0,0}=\frac{1}{2}(\beta_{1,-1}\alpha^{-1}+\beta_{-1,1}\alpha),
\end{equation}and
\begin{equation} \label{eq:1}
\zeta_0=\beta_{-1,1}\alpha.
\end{equation}
By (2.10) and $[\theta(L_n),\theta(L_m)] =(n-m)\theta(L_{m+n})\
(m+n\neq0),$ we can get
\begin{equation} \label{eq:1}
(n-m)\beta_{m+n,-(m+n)}
=n\beta_{n,-n}\alpha^{m}-m\alpha^{n}\beta_{m,-m} (m+n\neq0).
\end{equation} By (2.10) and $[\theta(L_n),\theta(L_{-n})] =2n\theta(L_{0})
+\frac{n-n^{3}}{12}\theta(C_L),$ we have
\begin{eqnarray*}
\frac{n-n^3}{12}\theta(C_L)&=&\frac{n-n^3}{12}C_L+[(n^2+n)\alpha^n\beta_{-n,n}-(n^2-n)\alpha^{-n}\beta_{n,-n}-2n\zeta_0]C_{LI}
\end{eqnarray*}\begin{equation} \label{eq:1}-(n\beta_{n,-n}\beta_{-n,n}+2n\gamma_0)C_I
+(n\alpha^n\beta_{-n,n}+n\alpha^{-n}\beta_{n,-n}-2n\beta_{0,0})I_0.\end{equation}
From (2.14) and (2.12), we can prove  by induction that
\begin{equation} \label{eq:1}
\beta_{m,-m}=\frac{m+1}{2}\alpha^{m-1}\beta_{1,-1}
-\frac{m-1}{2}\alpha^{m+1}\beta_{-1,1}.
\end{equation}
Setting $n=2$ in (2.15) and using (2.11), (2.13) and (2.16), we
get that
\begin{equation} \label{eq:1}
\theta(C_L)=C_L+12(\alpha^{-1}\beta_{1,-1}
-\alpha\beta_{-1,1})C_{LI}-(14\beta_{-1,1}\beta_{1,-1}-3\alpha^{-2}\beta^{2}_{1,-1}-3\alpha^{2}\beta^{2}_{-1,1})C_I.
\end{equation}
By Lemma 2.4 (ii)  and (2.9), we see  that $\lambda=1$ and
\begin{equation} \label{eq:1}
\theta(I_m)=\alpha_{m,-m}I_{-m}\ \mathrm{for}\ m\neq0,
\end{equation}
\begin{equation} \label{eq:1}
\theta(C_{I})=\alpha_{1,-1}\alpha_{-1,1} C_I.
\end{equation}
By (2.18) and  identity $\theta^2(I_m)=I_m,$ we have
\begin{equation} \label{eq:1}
\overline{\alpha_{m,-m}}\alpha_{-m,m}=1.
\end{equation}
Setting $n=\pm1$ in $[\theta(I_{-n}),\theta(L_{n})]
=-n\theta(I_{0})+(n^2-n)\theta(C_{LI})$ respectively and combing
with (2.10) and (2.18)-(2.19), we get that
\begin{equation} \label{eq:1}
\theta(I_0)=\alpha_{-1,1}\alpha I_0+2\alpha_{-1,1}\alpha
C_{LI}+\alpha_{-1,1}\beta_{1,-1}C_I,
\end{equation}
and
\begin{equation} \label{eq:1}
\theta(C_{LI})=-\alpha_{-1,1}\alpha
C_{LI}+\frac{1}{2}(\alpha_{1,-1}\alpha^{-1}-\alpha_{-1,1}\alpha)I_0+\frac{1}{2}(\alpha_{1,-1}\beta_{-1,1}
-\alpha_{-1,1}\beta_{1,-1})C_I.
\end{equation}
By $\theta^2(I_0)=I_0,$ we obtain that
\begin{equation} \label{eq:1}
\alpha\in\mathbb{R}^*,
\alpha\beta_{-1,1}+\overline{\beta_{1,-1}}\alpha_{-1,1}=0.
\end{equation}
By $[\theta(I_n),
\theta(L_{-n})]=n\theta(I_0)+(n^2+n)\theta(C_{LI})$, (2.16) and
(2.21)-(2.23), we can deduce that
\begin{equation} \label{eq:1}
\alpha_{n,-n}=\alpha_{-1,1}\alpha^{n+1}.
\end{equation}
Setting $m=1$ in (2.20) and combing with (2.24), we see that
$$\alpha_{1,-1}\overline{\alpha_{1,-1}}=\alpha^2.$$
Thus
\begin{equation} \label{eq:1}
\alpha_{1,-1}=\alpha e^{i\gamma},\
\alpha_{-1,1}=\alpha^{-1}e^{i\gamma}\ \mathrm{for}\ \mathrm{some}\
\gamma\in \mathbb{R}.
\end{equation}
From (2.23) and (2.25), we get that
\begin{equation} \label{eq:1}
\overline{\beta}_{1,-1}=-\alpha^2\beta_{-1,1} e^{-i\gamma},\
\overline{\beta}_{-1,1}=-\alpha^{-2}\beta_{1,-1} e^{-i\gamma}.
\end{equation}
By (2.17), (2.19), (2.22),(2.25)-(2.26) and identity
$\theta^2(C_L)=C_{L},$ we can deduce that
\begin{equation} \label{eq:1}
\beta_{1,-1}=0=\beta_{-1,1}.
\end{equation}
So type (i) follows from (2.9)-(2.27).

\vskip 3mm If $\bar{\theta}$ is of type (b), by a similar
discussion as that in type (i), we can prove that
\begin{equation} \label{eq:1}
\theta(L_n)
=-\alpha^{n}L_{n}+(\frac{n+1}{2}\alpha^{n-1}\beta_{1,-1}
-\frac{n-1}{2}\alpha^{n+1}\beta_{-1,1})I_{n}.
\end{equation}
\begin{equation} \label{eq:1}
\theta(C_L)=-C_L-12(\alpha^{-1}\beta_{1,-1}-\alpha\beta_{-1,1})C_{LI}
+(14\beta_{-1,1}\beta_{1,-1}-3\alpha^{-2}\beta^{2}_{1,-1}-3\alpha^{2}\beta^{2}_{-1,1})C_I.
\end{equation} where $\alpha\in S^{1}, \beta_{1,-1}, \beta_{-1,1}\in \mathbb{C}.$
\begin{equation} \label{eq:1}
\theta(I_n)=\alpha_nI_n \ \mathrm{for}\ n\neq0,
\end{equation}
where $\alpha_n=\alpha_{n,n}\in S^{1}(n\neq0).$
\begin{equation} \label{eq:1}
\theta(C_I)=-\alpha_1\alpha_{-1}C_I.
\end{equation}
\begin{equation} \label{eq:1}
\theta(I_0)=\alpha_{-1}\alpha I_0-\alpha_{-1}\beta_{1,-1}C_I.
\end{equation}
\begin{equation} \label{eq:1}
\theta(C_{LI})=\alpha_1\alpha^{-1}C_{LI}+\frac{1}{2}(\alpha_1\alpha^{-1}-\alpha_{-1}\alpha)I_0
+\frac{1}{2}(\alpha_{-1}\beta_{1,-1}-\alpha_1\beta_{-1,1})C_I.
\end{equation}
By $\theta^2(I_0)=I_0$ and $\theta^2(C_{LI})=C_{LI},$ we obtain
that
\begin{equation} \label{eq:1}
\alpha_1\overline{\beta_{1,-1}}=\overline{\alpha}\beta_{1,-1}, \
\alpha\beta_{-1,1}=\alpha_1\overline{\beta_{-1,1}}.
\end{equation}
By $[\theta(I_n),
\theta(L_{-n})]=n\theta(I_0)+(n^2+n)\theta(C_{LI})$, we can deduce
that
\begin{equation} \label{eq:1}
\alpha_{n}=\alpha_1\alpha^{n-1}.
\end{equation}
By (2.31)-(2.35) and identity $\theta^2(C_L)=C_L,$ we can deduce
that
\begin{equation} \label{eq:1}
(\overline{\alpha}-\alpha)\beta_{-1,1}=0,\
\beta_{1,-1}\beta_{-1,1}(1+\overline{\alpha}^2)=0.
\end{equation}
Thus type (ii) follows from (2.28)-(2.36). \hfill$\Box$

\vskip 3mm\noindent{\bf{Lemma 2.6.}}  Let $\theta$ be a
conjugate-linear anti-involution of  $\mathfrak{L}.$

(i) If $\theta=\theta_{\alpha,\gamma}^{+}.$  Then
$\theta(Vir)=Vir.$

(ii) If $\theta=\theta^-_{\alpha, \alpha_1,\beta_1,\beta_{-1}},$
we denote by $Vir^{'}$ the subalgebra of $\mathfrak{L}$ generated
by
$$\{C^{'}, L_{n}^{'}:=L_n+x_nI_{n}\mid
n\in\mathbb{Z}\},$$ where $x_n\in \mathbb{C}$ is determined by
$\overline{x_{n}}\alpha_1+\frac{n+1}{2}\beta_1-\frac{n-1}{2}\alpha^2\beta_{-1}=-\alpha
x_n,$ $C^{'}$ is determined by $\frac{n^3-n}{12}C^{'}=[L_n^{'},
L_{-n}^{'}]+2nL_0^{'}.$ Then $Vir^{'}\simeq Vir$ and
$\theta^-_{\alpha,
\alpha_1,\beta_1,\beta_{-1}}(L_{n}^{'})=-\alpha^{n}L_{n}^{'},
\theta^-_{\alpha, \alpha_1,\beta_1,\beta_{-1}}(C^{'})=-C^{'}.$

\vskip 3mm\noindent{\bf{Proof.}} It can be checked directly, we omit
the details.\hfill$\Box$

\vskip 3mm\noindent{\bf{Lemma 2.7.}} (Theorem 3.5 in Ref. [3]) Let
$V$ be a nontrivial irreducible weight $Vir$-module.

(i) If $V$ is unitary for some conjugate-linear anti-involution
$\theta$ of $Vir$,  then $\theta=\theta_{\alpha}^{+}$ for some
$\alpha>0.$

(ii) If $V$ is unitary for $\theta_{\alpha}^{+}$ for some
$\alpha>0,$ then $V$ is unitary for $\theta_{1}^{+}.$

\vskip 3mm\noindent{\bf{Proposition 2.8.}} Let $V$ be a nontrivial
irreducible weight $\mathfrak{L}$-module.

(i) If $V$ is unitary for some conjugate-linear anti-involution
$\theta$ of $\mathfrak{L}$, then
$\theta=\theta_{\alpha,\gamma}^{+}$ for some
$\alpha\in\mathbb{R}^+.$

(ii) If $V$ is unitary for $\theta_{\alpha,\gamma}^{+}$ for some
$\alpha>0,$ then $V$ is unitary for $\theta_{1,\gamma}^{+}.$

\vskip 3mm\noindent{\bf{Proof.}}   Suppose $V$ is unitary for some
conjugate-linear anti-involution $\theta$ of $\mathfrak{L}$. By
Lemma 2.6, $V$ can be viewed as a unitary $Vir^{'}$-module for the
conjugate-linear anti-involution $\theta|_{Vir^{'}}$, where
$Vir^{'}=Vir$ if $\theta=\theta^+_{\alpha,\gamma}$ and $Vir^{'}$
is that defined in Lemma 2.6(ii) if $\theta=\theta^-_{\alpha,
\alpha_1,\beta_1,\beta_{-1}}$. Then $V$ is a direct sum of
irreducible unitary $Vir^{'}$-modules since any unitary weight
$Vir^{'}$-module is complete reducible. We claim that $V$ is a
nontrivial $Vir^{'}$-module. Otherwise, for any $0\neq v\in V,$ by
$[L_1^{'},I_n]v=0 \ \mathrm{for}\ n\neq0,-1,$ we see that
$$I_nv=0, \forall n\neq0,1.$$ By $[L_2^{'},I_{-1}]v=0,$ we have
$$I_{1}v=0.$$ Since $I_{n}v=0 \ \mathrm{for} \ n\neq0,$ we have
$$C_Iv=0.$$
By $[L_{1}^{'},I_{-1}]v=0,$ we have $$I_0v=0.$$ By
$[L_{-1}^{'},I_1]v=0,$ we have $$C_{LI}=0.$$ So
$\mathfrak{L}.V=0,$ a contradiction. Thus there is a nontrivial
irreducible unitary $Vir^{'}$-submodule of $V$ for
conjugate-linear anti-involution $\theta|_{Vir^{'}}.$ By Lemma
2.7, $\theta|_{Vir^{'}}=\theta_{\alpha}^{+}$ for some $\alpha>0.$
Then by Proposition 2.5, we have
$\theta=\theta_{\alpha,\gamma}^{+}$ for some
$\alpha\in\mathbb{R}^+.$

The proof of (ii) is similar as that for Theorem 3.5 in Ref. [3],
we omit the details. \hfill$\Box$

\vskip 5mm\noindent{\section{Unitary representations for
$\mathfrak{L}$}}

In this section, we study the unitary  modules for $\mathfrak{L}$.
By Proposition 2.8, we see that the conjugate-linear
anti-involution is of the form $\theta_{\alpha, \gamma}^{+}$ for
some $\alpha\in\mathbb{R}^+.$  For the sake of simplicity, we
write $\theta_{\alpha, \gamma}^{+}$ by $\theta$. To prove the main
result, we give some lemmas first.

\vskip 2mm The following result is well known:

\noindent{\bf{Lemma 3.1.}} If $V$ is unitary weight module for
$Vir$, then $V$ is completely reducible.

\vskip 2mm In [12], all indecomposable Harish-Chandra modules of
intermediate series over the twisted Heisenberg-Virasoro algebra
are classified:

 \vskip 2mm\noindent{\bf{Lemma 3.2.}}(Theorem 3.5 in [12]) Let
 $V=\sum_{m\in\mathbb{Z}}\mathbb{F}v_m$ be an indecomposable
 $\mathfrak{L}$-module such that
 $\mathfrak{L}_mv_n\in\mathbb{C}v_{m+n}$ for all $m,n\in
 \mathbb{Z}.$ Then $V$ is isomorphic to one of the modules $\mathcal{A}_{a,b,c},
 \mathcal{A}(a,c), \mathcal{B}(a,c), \mathcal{U}_d, \mathcal{V}_d,\tilde{\mathcal{U}}_d,
 \tilde{\mathcal{V}}_c$ for appropriate $a,b,c,d\in\mathbb{F}.$

For more details, we refer the readers to Ref. [12].

It is well known that there are three types modules of the
intermediate series over $Vir$, denoted respectively by $A_{a,b},
A_{\alpha}, B_{\beta},$ they all have basis $\{v_k\mid
k\in\mathbb{Z}\}$ such that $C_L$ acts trivially and
\begin{align*}
&A_{a,b}: L_nv_k=(a+k+nb)v_{n+k};\\
&A_{\alpha}: L_nv_k=(n+k)v_{n+k}\  \mathrm{if} \ k\neq0,\ L_nv_0=n(n+a)v_n;\\
&B_{\beta}: L_nv_k=kv_{n+k}\  \mathrm{if} \ k\neq-n,\
L_nv_{-n}=-n(n+a)v_0.
\end{align*}for all $n, k\in \mathbb{Z}.$
About the  modules of the intermediate series of type $A_{a,b},$ we
have the following facts:

(1) $A_{a,b}$ is not simple if and only if $a\in \mathbb{Z}$ and
$b=0$ or $1.$

(2) $A_{a,b}\simeq A_{a^{'},b^{'}}$ if and only if (i) $a-a^{'}\in
\mathbb{Z}, b=b^{'}$ or (ii)$a-a^{'}\in \mathbb{Z}, a\notin
\mathbb{Z}, \{b,b^{'}\}=\{0,1\}.$

\vskip 3mm\noindent{\bf{Lemma 3.3.}} (Theorem 0.5 in [3]) Let $V$
be an irreducible unitary module of $Vir$ with finite-dimensional
weight spaces. Then either $V$ is highest or Lowest weight, or $V$
is isomorphic to $A_{a,b}$ for some $a\in \mathbb{R}, b\in
\frac{1}{2}+\sqrt{-1}\mathbb{R}.$

\vskip 3mm Now we prove the main theorem of this section.

\vskip 3mm\noindent{\bf{Theorem 3.4.}} Any unitary irreducible
Harish-Chandra module from intermediate series  with
conjugate-linear anti-involution $\theta_{\alpha,\gamma}^{+}$ over
twisted Heisenberg-Virasoro algebra $\mathfrak{L}$ is  a unitary
irreducible $\mathfrak{L}$-module of the intermediate series with
form $\mathcal{A}_{a,b,c}$ for $a\in \mathbb{R}, b\in
\frac{1}{2}+\sqrt{-1}\mathbb{R}, c\in \mathbb{C}$ satisfying
$c=\overline{c}e^{-i\gamma}$.

\vskip 2mm\noindent{\bf{Proof.}} Let $V$ be a unitary irreducible
$\mathfrak{L}$-module of intermediate series for a
conjugate-linear anti-involution $\theta.$ Then the central
elements $C_L, C_I, C_{LI},I_0$ are assigned zero by [12]. By
Lemma 2.6(i) and Proposition 2.8, $V$ is also unitary for $Vir$
and then by Lemma 3.1 and Lemma 3.3, we can suppose that
$$V=A_{a_1,b_1}\oplus\cdots \oplus A_{a_K,b_K}\oplus W,$$
where $a_i\in\mathbb{R}, b_i\in\frac{1}{2}+\sqrt{-1}\mathbb{R},$
$W$ is a trivial $Vir$-module.

\vskip 3mm\noindent{\bf{Case 1.}} $W=0.$

\vskip 3mm\noindent In this subcase
$$V=A_{a_1,b_1}\oplus\cdots \oplus A_{a_K,b_K}.$$
Let $\{v_k\mid k\in \mathbb{Z}\}$ be a basis of $A_{a_1,b_1}$ such
that $L_nv_k=(a_1+k+nb_1)v_{n+k}.$ As a irreducible
$\mathfrak{L}$-module, $V$ is generated by the $L_0$-eigenvector
$v_0$ with eigenvalue $a_1.$ Thus $L_0$-eigenvalue on $V$ are of the
form $a_1+n, n\in \mathbb{Z}.$ This means that $a_i\in \{a_1+n\mid
n\in \mathbb{Z}\},$ $i=1,\cdots,K.$ Recall that $A_{a+n,b}\simeq
A_{a,b}$ for any $n\in \mathbb{Z}.$  So there exists $0\leq a<1$
such that $A_{a_i,b_i}$ are of the form $A_{a,b_i}.$ i.e.,
\begin{equation} \label{eq:1}
V=A_{a,b_1}\oplus\cdots \oplus A_{a,b_K},
\end{equation}
where $0\leq a<1, b_i\in\frac{1}{2}+\sqrt{-1}\mathbb{R}.$

\vskip 3mm\noindent{\bf{\emph{Claim 1.}}} $V=A_{a,b}.$

\vskip 3mm\noindent{\bf{\emph{proof of claim 1}.}} By (3.1) we can
choose a basis of $V:$
$$\{v_{k,l}\mid k\in\mathbb{Z}, 1\leq l\leq K\}$$
such that
\begin{equation} \label{eq:1}
L_{m}v_{k,l}=(a+k+mb_l)v_{k+m,l},
\end{equation}
for $m\in\mathbb{Z}, 1\leq l\leq K.$ Suppose \begin{equation}
\label{eq:1}
I_{1}v_{k,l}=\sum_{l^{'}=1}^{K}\mu_{k,l}^{l^{'}}v_{k+1,l^{'}}.
\end{equation}
By (3.2), (3.3) and $[L_{-1}^{'},I_1]=I_0+2C_{LI},$ we have
$$
I_0v_{k,l}=\sum_{l^{'}=1}^{K}[(a+k+1-b_{l^{'}})\mu_{k,l}^{l^{'}}-(a+k-b_{l})\mu_{k-1,l}^{l^{'}}]v_{k,l^{'}}.
$$
Considering $I_0V=cV$ for some $c\in\mathbb{C}$ since
$I_0\in\mathfrak{C},$ we see that
\begin{equation}\label{eq:1}
(a+k+1-b_{l^{'}})\mu_{k,l}^{l^{'}}=(a+k-b_{l})\mu_{k-1,l}^{l^{'}}\
\mathrm{for}\ l^{'}\neq l,
\end{equation}
and
\begin{equation} \label{eq:1}
I_0v_{k,l}=[(a+k+1-b_{l})\mu_{k,l}^{l}-(a+k-b_{l})\mu_{k-1,l}^{l}]v_{k,l}=cv_{k,l}.
\end{equation}

If there exists $k_0$ such that $a+k_0-b_l=0$ or
$a+k_0+1-b_{l^{'}}=0,$ noting that $0\leq a<1, b_l, b_{l^{'}}\in
\frac{1}{2}+\sqrt{-1}\mathbb{R},$ we can deduce recursively from
(3.4) that $\mu_{k,l}^{l^{'}}=0, \forall k\in\mathbb{Z}, l\neq
l^{'}.$

If $a+k-b_l\neq0$ and $a+k+1-b_{l^{'}}\neq0$ for all
$k\in\mathbb{Z},$ then for $l\neq l^{'},$
\begin{equation}\label{eq:1}
\mu_{k,l}^{l^{'}}=\frac{a+k-b_{l}}{a+k+1-b_{l^{'}}}\mu_{k-1,l}^{l^{'}},\
\mu_{k-1,l}^{l^{'}}=\frac{a+k+1-b_{l^{'}}}{a+k-b_{l}}\mu_{k,l}^{l^{'}}.
\end{equation}
By (3.2), (3.3) and $[L_{-2},I_1]=I_{-1}+6C_{LI},$ we have
\begin{equation} \label{eq:1}
I_{-1}v_{k,l}=\sum_{l^{'}=1}^{K}[(a+k+1-2b_{l^{'}})\mu_{k,l}^{l^{'}}-(a+k-2b_{l})\mu_{k-2,l}^{l^{'}}]v_{k-1,l^{'}}.
\end{equation}
By (3.2), (3.7) and  $[L_1,I_{-1}]=-I_0,$ we have
\begin{eqnarray*}
I_0v_{k,l}&=&\sum_{l^{'}=1}^{K}\{(a+k+b_l)[(a+k+2-2b_{l^{'}})\mu_{k+1,l}^{l^{'}}-(a+k+1-2b_{l})\mu_{k-1,l}^{l^{'}}]
\end{eqnarray*}
\begin{equation} \label{eq:1}
-(a+k-1+b_{l^{'}})[(a+k+1-2b_{l^{'}})\mu_{k,l}^{l^{'}}-(a+k-2b_{l})\mu_{k-2,l}^{l^{'}}]\}v_{k,l^{'}}.
\end{equation}
Comparing (3.5) and (3.8), we have
\begin{eqnarray*}
(a+k+b_l)[(a+k+2-2b_{l^{'}})\mu_{k+1,l}^{l^{'}}-(a+k+1-2b_{l})\mu_{k-1,l}^{l^{'}}]\\
\end{eqnarray*}
\begin{equation} \label{eq:1}
-(a+k-1+b_{l^{'}})[(a+k+1-2b_{l^{'}})\mu_{k,l}^{l^{'}}-(a+k-2b_{l})\mu_{k-2,l}^{l^{'}}]=0\
\mathrm{for}\ l\neq l^{'}.
\end{equation}
By (3.6) and (3.9), we obtain that
\begin{eqnarray*}
&&[(a+k+b_{l})(a+k-b_{l})(a+k+2-2b_{l^{'}})(a+k+1-b_l)(a+k-1-b_l)\\
&&-(a+k+b_{l})(a+k+1-2b_l)(a+k+1-b_{l^{'}})(a+k+2-b_{l^{'}})(a+k-1-b_l)\\
&&-(a+k-1+b_{l^{'}})(a+k+1-2b_{l^{'}})(a+k+2-b_{l^{'}})(a+k-b_l)(a+k-1-b_l)\\
&&+(a+k-1+b_{l^{'}})(a+k-2b_l)(a+k-b_{l^{'}})(a+k+1-b_{l^{'}})(a+k+2-b_{l^{'}})]\mu_{k,l}^{l^{'}}\\
&&=0.
\end{eqnarray*}
Computing the above equality by mathematical software, such as
maple, we get that:
\begin{equation} \label{eq:1}
\{[-(b_l+2-b_{l^{'}})(b_l+1-b_{l^{'}})((b_l+b_{l^{'}})(b_l+b_{l^{'}}-1)-2b_{l^{'}})]k+f(b_l,b_{l^{'}})\}\mu_{k,l}^{l^{'}}=0,
\end{equation}
where $f(b_l,b_{l^{'}})\in\mathbb{C}$ is a constant determined by
$b_l$ and $b_{l^{'}}.$ Note that $b_{l},
b_{l^{'}}\in\frac{1}{2}+\sqrt{-1}\mathbb{R},$ it is easy to deduce
that the coefficient of $k$ in (3.10) is nonzero. Thus there
exists $k\in \mathbb{Z}$ such that $\mu_{k,l}^{l^{'}}=0.$ Then, by
(3.6), we have
\begin{equation} \label{eq:1}
\mu_{k,l}^{l^{'}}=0,\forall k\in\mathbb{Z}, l\neq l^{'}.
\end{equation}
From (3.2), (3.3) and (3.11), we see that each subspace
$A_{a,b_l}(1\leq l\leq K)$ in (3.1) is a submodule of $V.$ Since $V$
is irreducible, we have $ V=A_{a,b}. $ Thus Claim 1 holds.

\vskip3mm By Claim 1, for the sake of simplicity, we can rewrite
(3.2), (3.3), (3.5) and (3.8) as following:
\begin{equation} \label{eq:1}
L_{m}v_{k}=(a+k+mb)v_{k+m},
\end{equation}
\begin{equation} \label{eq:1}
I_{1}v_{k}=\mu_{k}v_{k+1}.
\end{equation}
\begin{equation} \label{eq:1}
I_0v_{k}=[(a+k+1-b)\mu_{k}-(a+k-b)\mu_{k-1}]v_{k}=cv_{k}.
\end{equation}
\begin{eqnarray*}
I_0v_{k}=\{(a+k+b)[(a+k+2-2b)\mu_{k+1}-(a+k+1-2b)\mu_{k-1}]
\end{eqnarray*}
\begin{equation} \label{eq:1}
-(a+k-1+b)[(a+k+1-2b)\mu_{k}-(a+k-2b)\mu_{k-2}]\}v_{k}=cv_k.
\end{equation}

\vskip 3mm\noindent{\bf{\emph{Claim 2.}}} $I_mv_k=cv_{m+k}, \forall
m,k\in\mathbb{Z}.$

\vskip 3mm\noindent{\bf{\emph{proof of claim 2}.}} If $(a+k+1-b)=0,$
we see that $a=\frac{1}{2}=b, k=-1$ since $0\leq a<1,
b\in\frac{1}{2}+\sqrt{-1}\mathbb{R}.$ Thus by (3.14) we have
$\mu_{-2}=c.$ Then we can deduce recursively that $\mu_k=c, \forall
k\in\mathbb{Z}.$ Similarly, if $(a+k-b)=0,$ we can also deduce that
$\mu_k=c, \forall k\in\mathbb{Z}.$

Now we suppose $(a+k+1-b)\neq0$ and $(a+k-b)\neq0$ for all
$k\in\mathbb{Z}.$ By (3.14) we have
\begin{equation} \label{eq:1}
\mu_k=\frac{c+(a+k-b)\mu_{k-1}}{a+k+1-b},
\mu_{k-1}=\frac{(a+k+1-b)\mu_{k}-c}{a+k-b}.
\end{equation}
By (3.15) and (3.16) we get that
\begin{eqnarray*}
&&[(a+k-b)(a+k-1-b)(a+k+b)(a+k+2-2b)(a+k+1-b)\\
&&-(a+k+2-b)(a+k-1-b)(a+k+b)(a+k+1-2b)(a+k+1-b)\\
&&-(a+k-1+b)(a+k+1-2b)(a+k+2-b)(a+k-b)(a+k-1-b)\\
&&+(a+k+2-b)(a+k-b)(a+k-1+b)(a+k-2b)(a+k+1-b)]\mu_k\\
&&=(a+k+2-b)(a+k-b)(a+k-1-b)c\\
&&-(a+k-b)(a+k-1-b)(a+k+2-2b)(a+k+b)c\\
&&-(a+k+2-b)(a+k-1-b)(a+k+b)(a+k+1-2b)c\\
&&+2(a+k+2-b)(a+k-b)(a+k-1+b)(a+k-2b)c.\\
\end{eqnarray*}
Resort to maple software again, we have
$$
[(8b-8b^2)k+4b+8ab-12b^2+8b^3-8ab^2](\mu_k-c)=0.
$$
Noting that $b\in\frac{1}{2}+\sqrt{-1}\mathbb{R}$ and combing with
(3.16), we can easily deduce that
\begin{equation} \label{eq:1}
\mu_k=c,\ \forall k\in\mathbb{Z}.
\end{equation}
By (3.12),(3.13) and (3.17), we have $
I_{m+1}v_k=[L_m,I_1]v_k=cv_{m+k+1},\ \forall m,k\in\mathbb{Z}. $
Thus Claim 2 holds.

From Claim 1, Claim 2 and (3.12) we see that $V$ is a unitary
irreducible $\mathfrak{L}$-module of the intermediate series with
form $\mathcal{A}_{a,b,c}$ in  Case 1. Moreover, by
$$\langle I_0v_k,I_0v_k\rangle=c\overline{c}\langle
v_k,v_k\rangle$$ and $$\langle I_0v_k,I_0v_k\rangle=\langle
v_k,e^{i\gamma}c^2v_k\rangle=\overline{c}^2e^{-i\gamma}\langle
v_k,v_k\rangle,$$ we see that $$c=\overline{c}e^{-i\gamma}.$$

\vskip 3mm\noindent{\bf{Case 2.}} $W\neq0.$

\vskip 3mm\noindent Choose an arbitrary nonzero element $w\in W.$ We
see that $V$ is generated by $w$ since  $V$ is an irreducible
$\mathfrak{L}$-module. If $I_1.w=0,$ then $\mathfrak{L}.W=0$ and $V$
is a trivial $\mathfrak{L}$-module, a contradiction.  Thus
$I_1w\neq0.$ By $L_0I_1w=I_1w,$ we see that
$$
I_1w\in A_{a_1,b_1}\oplus\cdots \oplus A_{a_K,b_K}.
$$
Moreover,
$$
A_{a_i,b_i}\simeq A_{0,b_i}, \forall i\in\{1,\cdots, K\}
$$
since $V$ is generated by the eigenvector $w$ of $L_0$ with
eigenvalue $0.$ Thus
$$V=A_{0,b_1}\oplus\cdots\oplus A_{0,b_K}\oplus W.$$
Choose the standard basis $\{v_{k,i}\mid k\in\mathbb{Z}\}$ for each
$A_{0,b_i}.$ Suppose
$$
I_1v_{k,l}=\sum_{l^{'}=1}^{K}\mu_{k,l}^{l^{'}}v_{k+1,l^{'}}+w_{k,l},
$$
where $w_{k,l}\in W.$ By a similar calculation as that of Claim 1
in Case 1, we can deduce that
$$V=A_{0,b}.$$
This contradicts with the assumption  that $W\neq 0$. Thus Case 2
is impossible.

This completes the proof of Theorem 3.4. \hfill$\Box$

\vskip 3mm\noindent{\bf{Acknowledgments}}

The research was done during the visit of the first author to
Chern Institute of Mathematics in 2011. The hospitality and
financial support of Chern Institute of Mathematics are gratefully
acknowledged. It is also a pleasure to thank Professor Chengming
Bai for useful conversations.

\end{document}